\documentclass[12pt]{article}
\usepackage{amstex}
\usepackage{amssymb}
\usepackage{amscd}
\usepackage{a4}
\usepackage{pb-diagram}

\def\fc{\displaystyle{\buildrel \subset\over \longrightarrow}}
\def\Talpha#1{\vbox{\ialign{##\crcr
$\alpha$\crcr\noalign{\kern2pt\nointerlineskip}
$\hfil\displaystyle{#1}\hfil$\crcr}}} 
 
\def\Onabla#1{\vbox{\ialign{##\crcr$\,\scriptstyle{0}$\crcr
\noalign{\kern2pt\nointerlineskip}
$\hfil\displaystyle{#1}\hfil$\crcr}}}
\def\im{\mbox{Im}}

 \def\cale{{\cal E}}
 \def\cala{{\cal A}}
\def\fraca{{\mathfrak{A}}}

\def\calm{{\cal M}}

\def\cals{{\cal S}}
\def\calh{{\cal H}}
\def\fm{\mathfrak m}
\def\ft{\mathfrak T}
 \def\oplusinf{\mathop{\oplus}}
\def\otimesinf{\mathop{\otimes}} \def\cals{{\cal S}}

\def\cals{{\cal S}}

\def\bbbone{\mbox{\rm 1\hspace {-.6em} l}}

\def\zn{\mathbb{Z}_N}
\def\fil{\displaystyle{\buildrel {[i^\ell]}\over \longrightarrow}}
\def\fdm{\displaystyle{\buildrel {[d^m]}\over \longrightarrow}}

 \newtheorem{theorem}{THEOREM}
\newtheorem{lemma}{LEMMA} \newtheorem{definition}{DEFINITION}
\newtheorem{proposition}{PROPOSITION} 
\begin{document}

\baselineskip=0.9cm \begin{center} {\Large\bf UNIVERSAL $q$-DIFFERENTIAL
CALCULUS
AND $q$-ANALOG OF HOMOLOGICAL ALGEBRA} \end{center} \vspace{0.75cm}

\begin{center} Michel DUBOIS-VIOLETTE \\ \vspace{0.3cm} {\small Laboratoire de
Physique Th\'eorique et Hautes Energies\footnote{Laboratoire associ\'e au
Centre
National de la Recherche Scientifique - URA D0063}\\ Universit\'e Paris XI,
B\^atiment 211\\ 91 405 Orsay Cedex, France\\ flad$@$qcd.th.u-psud.fr}\\ and\\
Richard KERNER\\ \vspace{0.3cm} {\small Laboratoire Physique Th\'eorique,
Gravitation et Cosmologie Relativistes \footnote{Laboratoire associ\'e au
Centre
National de la Recherche Scientifique - URA D0769}\\ Tour 22/12 - 4e \'etage -
Bo\^\i te 142, 4 place Jussieu\\ 75005 Paris, France\\ rk$@$ccr.jussieu.fr}

\end{center} \vspace{1cm}

\begin{center} \today \end{center}

\vspace {1cm}

\noindent L.P.T.H.E.-ORSAY 96/48\\

\newpage \begin{abstract}

We recall the definition of $q$-differential algebras and discuss some
representative examples. In particular we construct the $q$-analog of the
Hochschild coboundary.  We then construct the universal $q$-differential
envelope
of a unital associative algebra and study its properties. The paper also
contains
general results on $d^N=0$.

\end{abstract}

\section{Introduction and algebraic preliminaries}

At the origin of this paper there is the long-standing physically-motivated
interest of one of the authors (R.K.) on ${\mathbb Z}_3$-graded structures and
differential calculi \cite{RK} although here the point of view is somehow
different. There is also the observation that the simplicial (co)-homology
admits
${\mathbb Z}_N$ versions leading to cyclotomic homology \cite{Sark} and that,
more generally, this suggests that one can introduce ``$q$-analog of
homological
algebra" for each primitive root $q$ of the unity \cite{Kapr}. Moreover the
occurrence of various notions of ``$q$-analog" in connection with quantum
groups
suggests to include in the formulation the general case where $q$ is not
necessarily a root of unity but is an arbitrary invertible complex number
\cite{Kapr}. It is our aim here to go further in this direction.\\

Throughout this paper, we shall be interested in complex associative graded
algebras equipped with endomorphisms $d$ of degree 1 satisfying a twisted
Leibniz
rule, {\sl the $q$-Leibniz rule}, of the form
$d(\alpha\beta)=d(\alpha)\beta+q^{\partial\alpha}\alpha d(\beta)$, where $q$ is
a
given complex number distinct of 0 and where $\partial\alpha$ denotes the
degree
of $\alpha$. Furthermore, whenever $q^N=1$, for an integer $N\geq 1$, we shall
add the rule $d^N=0$. We shall refer to $d$ as {\sl the $q$-differential} of
the
graded algebra. Thus an ordinary differential on a graded algebra is just a
$(-1)$-differential in this terminology. Our aim is to produce universal
objects
in this class. Before entering the subject we want to discuss two problems
connected with the case $q^N=1$, i.e. the case where $q$ is a primitive root of
the unity.\\

We shall be concerned here with the case of $\mathbb{N}$-graded algebras.
However
when $q^N=1$, it is very natural to consider graduation over
$\mathbb{Z}_N=\mathbb{Z}/N\mathbb{Z}$ instead of over $\mathbb{N}$. So let us
recall how one can identify a $\mathbb{Z}_N$-graded algebra with a
$\mathbb{N}$-graded one. Let $\fraca=\oplusinf_{p\in\zn}\fraca^p$ be a
$\zn$-graded algebra and let $n\mapsto p(n)$ be the canonical projection of
$\mathbb{N}$ onto $\zn$. We associate to $\fraca$ a $\mathbb{N}$-graded algebra
$p^\ast\fraca=\displaystyle{\oplusinf_{n\in \mathbb{N}}}p^\ast\fraca^n$ in the
following manner. A homogeneous element of $p^\ast\fraca$ is a pair
$(n,\alpha)\in \mathbb{N}\times \fraca$ of an integer $n\in\mathbb{N}$ and an
homogeneous element $\alpha$ of $\fraca$ such that $\partial\alpha=p(n)$ and we
identify $p^\ast\fraca^n=(n,\fraca^{p(n)})$ with the vector space
$\fraca^{p(n)}$. The product in $p^\ast\fraca=\oplus p^\ast\fraca^n$ is then
defined by $(m,\alpha)(n,\beta)=(m+n,\alpha\beta)$. The canonical projection
$\pi:p^\ast\fraca\rightarrow \fraca$ defined by $\pi(n,\alpha)=\alpha$ is an
algebra homomorphism which is graded in the sense that one has
$\pi(p^\ast\fraca^n)\subset \fraca^{p(n)}$. The $\mathbb{N}$-graded algebra
$p^\ast\fraca$ is characterized, up to an isomorphism, by the following
universal
property: Any graded homomorphism $\varphi$ of a $\mathbb{N}$-graded algebra
$\Omega$ into $\fraca$ factorizes through a unique homomorphism of
$\mathbb{N}$-graded algebra $\bar\varphi:\Omega\rightarrow p^\ast\fraca$ as
$\varphi=\pi\circ \bar\varphi$. Let $D$ be a homogeneous linear mapping of
$\fraca$ into itself and let $k$ be the unique positive integer strictly
smaller
than $N$ such that $p(k)$ is the degree of $D$. Then there is a unique linear
mapping $p^\ast(D)$ of $p^\ast\fraca$ into itself which is homogeneous of
degree
$k$ and satisfies $\pi\circ p^\ast(D)=D\circ\pi$. The construction of
$p^\ast(D)$
is obvious.\\

As already stressed, we impose $d^N=0$ whenever $q^N=1$. More generally let $E$
be a vector space equipped with an endomorphism $d$ satisfying $d^N=0$, $N$
being
an integer greater than or equal to 2. For each integer $k$ with $0\leq k\leq
N$,
one has $\im(d^{N-k})\subset \ker(d^k)$ so the vector space
$H^{(k)}=\ker(d^k)/\im(d^{N-k})$ is well defined. One has $H^{(0)}=H^{(N)}=0$
and
the $H^{(k)}$ for $1\leq k\leq N-1$ are the generalized homologies of $E$. Let
$\ell$ and $m$ be two positive integers such that $\ell+m\leq N$. The inclusion
$i^\ell$ : $\ker(d^m)\subset \ker(d^{\ell+m})$ induces a linear mapping
$[i^\ell]:H^{(m)}\rightarrow H^{(\ell+m)}$ since $\im(d^{N-m})\subset
\im(d^{N-(\ell+m)})$. On the other hand, one has
$d^m(\ker(d^{\ell+m}))\subset\ker(d^\ell)$ and $d^m(\im(d^{N-(\ell+m)}))\subset
\im(d^{N-\ell})$ and therefore $d^m$ induces a linear mapping
$[d^m]:H^{(\ell+m)}\rightarrow H^{(\ell)}$. One has the following result.

\begin{lemma}

The hexagon $(\calh^{\ell,m})$ of homomorphisms
\[ \begin{diagram} \node{}
\node{H^{(\ell+m)}} \arrow{e,t}{[d^m]} \node{H^{(\ell)}}
\arrow{se,t}{[i^{N-(\ell+m)}]} \node{} \\ \node{H^{(m)}} \arrow{ne,t}{[i^\ell]}
\node{} \node{} \node{H^{(N-m)}} \arrow{sw,b}{[d^\ell]} \\ \node{}
\node[1]{H^{(N-\ell)}} \arrow{nw,b}{[d^{N-(\ell+m)}]} \node{H^{(N-(\ell+m))}}
\arrow{w,b}{[i^m]} \node{} \end{diagram} \]
is exact.
\end{lemma}
\noindent
\underbar{Proof}. It is clearly sufficient to show that the sequences
$H^{(m)}\fil H^{(\ell+m)}\fdm H^{(\ell)}$ and $H^{(\ell+m)} \fdm H^{(\ell)}
\displaystyle{\buildrel {[i^{(N-(\ell+m))}]}\over \longrightarrow}H^{(N-m)}$
are
exact. It is straightforward that $[d^m]\circ [i^\ell]$ is the zero mapping of
$H^{(m)}$ into $H^{(\ell)}$ and that $[i^{N-(\ell+m)}]\circ [d^m]$ is the zero
mapping of $H^{(\ell+m)}$ into $H^{(N-m)}$. Let $c\in \ker(d^{(\ell+m)})$ be
such
that $d^m c=d^{N-\ell}c'$ for some $c'\in E$; Then $d^m(c-d^{N-(\ell+m)}c')=0$
i.e. $[c]\in H^{(\ell+m)}$ is in $[i^\ell](H^{(m)})$ which achieves the proof
of
the exactness of the first sequence. Let $c\in \ker(d^\ell)$ be such that
$c=d^mc'$ for some $c'\in E$; then one has $d^{\ell+m}c'=0$ which means that
$[c]\in H^{(\ell)}$ is in $[d^m](H^{(\ell+m)})$ which achieves the proof of the
exactness of the second sequence.$\square$\\ Notice that the content of this
lemma is nontrivial only if $\ell\geq 1,\ m\geq 1$ and $N-(\ell+m)\geq 1$ which
implies $N\geq 3$. In the case $N=3$ the only nontrivial choice is
$(\ell,m)=(1,1)$ so, in this case there is only one (nontrivial) hexagon,
namely
$(\calh^{1,1})$. In the case $N=4$, there are 3 possible choices for $(\ell,m)$
namely $(\ell,m)=(1,1),(\ell,m)=(1,2)$ and $(\ell,m)=(2,1)$. However it is
readily seen that, (for $N=4$), $(\calh^{1,1}),(\calh^{1,2})$ and
$(\calh^{2,1})$
are identical; one passes from one to the others by applying ``rotations of
$2\pi/3$". Thus, for a given integer $N\geq 3$, it is not completely obvious to
count the number of independent nontrivial hexagons. In any case, this lemma is
very useful for the computations. Practically we shall apply it in the graded
case where $E=\oplus_{n\in \mathbb Z}E^n$ is a $\mathbb Z$-graded vector space
and where $d$ is homogeneous of degree 1, (i.e. $d(E^n)\subset E^{n+1}$). In
this
case, the hexagon $(\calh^{\ell,m})$ of the lemma splits into $N$ long exact
sequences $(\cals^{\ell,m}_p),\ p\in\{0,1,\dots,N-1\}$. \[ \begin{array}{ll} &
\cdots\ \longrightarrow  H^{(m),Nr+p}\ \displaystyle{\buildrel {[i^\ell]}\over
\longrightarrow}\ H^{(\ell+m),Nr+p}\ \displaystyle{\buildrel {[d^m]}\over
\longrightarrow}\ H^{(\ell),Nr+p+m\ }\displaystyle{\buildrel
{[i^{N-(\ell+m)}]}\over \longrightarrow}\\ \\ & H^{(N-m),Nr+p+m}\
\displaystyle{\buildrel {[d^\ell]}\over \longrightarrow}\
H^{(N-(\ell+m)),Nr+p+\ell+m}\ \displaystyle{\buildrel {[i^m]}\over
\longrightarrow}\ H^{(N-\ell),Nr+p+\ell+m}\\ \\ & \displaystyle{\buildrel
{[d^{N-(\ell+m)}]}\over \longrightarrow}\ H^{(m),N(r+1)+p}\
\displaystyle{\buildrel {[i^\ell]}\over \longrightarrow}\ \cdots \end{array} \]
\begin{flushright} $(\cals^{\ell,m}_p)$ \end{flushright} where
$H^{(k),n}=\{x\in
E^n\vert d^k(x)=0\}/d^{N-k}(E^{n+k-N})$. Notice that, if instead of being
graded
over $\mathbb Z$, $E$ is graded over $\mathbb Z_N$ then the $N$ exact sequences
$(\cals^{\ell,m}_p)$ are again $N$ exact hexagons because in $\mathbb Z_N$ one
has $Nr+p=N(r+1)+p=p$.\\

In degree 0, the $q$-Leibniz rule reduces to the ordinary Leibniz rule. Thus a
$q$-differential induces a derivation of the subalgebra of elements of degree 0
into the space of elements of degree 1 which is a bimodule over the algebra of
elements of degree 0. In this context, let us recall the construction of the
universal derivation \cite{CE}. Let $\cala$ be a unital associative algebra and
let $\Omega^1(\cala)$ be the kernel of the product
$\fm:\cala\otimes\cala\rightarrow \cala$ of $\cala$, $\fm(x\otimes y)=xy$. The
mapping $\mathfrak m$ is a bimodule homomorphism so $\Omega^1(\cala)$ is a
bimodule over $\cala$. One defines a derivation $d$ of $\cala$ into
$\Omega^1(\cala)$ by $d(x)=\bbbone \otimes x-x\otimes \bbbone$ for $x\in
\cala$.
The derivation $d$ is universal in the sense that {\sl for any derivation $X$
of
$\cala$ into a bimodule $\calm$ over $\cala$, there is a unique bimodule
homomorphism $i_X$ of $\Omega^1(\cala)$ into $\calm$ such that $X=i_X\circ d$}.
This universal property characterizes the pair $(\Omega^1(\cala),d)$ uniquely,
up
to an isomorphism. We proceed now to recall the construction of the universal
differential calculus over $\cala$ \cite{Kar}. Set $\Omega^0(\cala)=\cala$ and
$\Omega^n(\cala)=\otimes^n_{\cala}\Omega^1(\cala)$. The direct sum
$\Omega(\cala)=\oplus_n\Omega^n(\cala)$ is an associative graded algebra for
the
tensor product over $\cala$; it is in fact the tensor algebra over $\cala$ of
the
bimodule $\Omega^1(\cala)$. The derivation\linebreak[4] $d:\cala\rightarrow
\Omega^1(\cala)$ extends uniquely into a differential (i.e. a
$(-1)$-differential) of $\Omega(\cala)$ which will be again denoted by $d$.
Thus,
$\Omega(\cala)$ is a graded differential algebra, (i.e. a graded
$(-1)$-differential algebra in the sense of the definition of  next section).
This graded differential algebra is characterized, up to an isomorphism, by the
following universal property~: {\sl Any homomorphism of associative unital
algebra $\varphi$ of $\cala$ into the algebra $\fraca^0$ of the elements of
degree 0 of a graded differential algebra $\fraca=\oplus_{n\in \mathbb
N}\fraca^n$ extends uniquely into a homomorphism of graded differential algebra
$\bar\varphi:\Omega(\cala)\rightarrow \fraca$}. This is why $\Omega(\cala)$ is
called {\sl the universal differential envelope} of $\cala$ or {\sl the
universal
differential calculus over $\cala$}. It is one of the aims of this paper to
generalize this construction (corresponding to $q=-1$) for the $q$-differential
calculus.

\section{$q$-differential calculus}

In the rest of the paper $q$ is a complex number with $q\not=0$ and we shall
use
the following definition.\\

\begin{definition}{\sl A graded $q$-differential algebra} is a $\mathbb
N$-graded
unital $\mathbb C$-algebra $\fraca=\displaystyle{\oplusinf_{n\in \mathbb N}}$
$\fraca^n$ equipped with an endomorphism $d$ of degree one satisfying
$d(\alpha\beta)=d(\alpha)\beta+q^a\alpha d(\beta)$, $\forall \alpha\in
\fraca^a$
and $\forall \beta\in \fraca$, and such that $d^N=0$ whenever $q^N=1$ for
$N\in\mathbb N$ with $N\not= 0$. Let $\cala$ be a unital $\mathbb C$-algebra.
{\sl A $q$-differential calculus  over $\cala$} is a graded $q$-differential
algebra $\fraca=\displaystyle{\oplusinf_n}\fraca^n$ such that $\cala$ is a
subalgebra of $\fraca^0$. \end{definition}

Notice that a graded 1-differential algebra is just a $\mathbb N$-graded
algebra
$(d=0)$, that a graded ($-1$)-differential algebra is just a  $\mathbb
N$-graded
differential algebra in the usual sense and that, if $\fraca=\oplus\fraca^n$ is
a
$q$-differential calculus over $\cala$ with $q\not=1$, then the restriction of
$d$ to $\cala$ is just a derivation of $\cala$ into the bimodule $\fraca^1$
over
$\cala$.

Let us introduce, as usual, the $q$-analogs of basic numbers, of factorials and
of binomial coefficients \[ \begin{array}{ll} [n]_q & =
\frac{1-q^n}{1-q}=1+q+\dots + q^{n-1},\\ \\ \mbox{[}n!]_q & = [1]_q[2]_q\dots
[n]_q\ \mbox{and}\\ \\ \left[\begin{array}{l} n\\p\end{array}\right]_q & =
\frac{[n!]_q}{[p!]_q[(n-p)!]_q} \end{array} \] where $n,p\in \mathbb N$ and
$n\geq p$. By induction on $n$, it follows from the $q$-Leibniz rule
$d(\alpha\beta)=d(\alpha)\beta+q^a\alpha d(\beta)$ that one has~:
\begin{equation} d^n(\alpha\beta)=\sum^n_{p=0}
q^{ap}\left[\begin{array}{l}n\\p\end{array}\right]_q d^{n-p}(\alpha)d^p(\beta)
\label{eq1} \end{equation} for $\alpha\in\fraca^a$ and $\beta\in \fraca$,
($n\in
\mathbb N$). It is worth noticing here that the consistency of $d^N=0$ whenever
$q^N=1$ with the $q$-Leibniz rule follows from the fact that (\ref{eq1})
implies
for $q^N=1$ that one has $d^N(\alpha\beta)=d^N(\alpha)\beta+\alpha
d^N(\beta)$.\\

There is an obvious notion of homomorphism of $q$-differential algebra. Given a
unital algebra $\cala$, a morphism of a $q$-differential calculus over $\cala$
into another one is a homomorphism of the corresponding $q$-differential
algebra
which induces the identity mapping of $\cala$ onto itself. It is the aim of
Section 4 to produce an initial object in the category of $q$-differential
calculi over $\cala$, (i.e. a universal graded $q$-differential envelope for
$\cala$). In the remaining part of this section, we present some examples.\\

\subsection*{Example 1: Matrix algebra $M_N(\mathbb C)$}

Let $N\in \mathbb N$ with $N\geq 2$  and let $q$ be a primitive $N$-root of the
unity, (e.g. $q=\exp\left(\frac{2\pi i}{N}\right)$). Let us introduce the usual
standard basis $E^k_\ell,\ (k,\ell\in \{1,\dots,N\})$, of the matrix algebra
$M_N(\mathbb C)$ defined by $(E^k_\ell)^i_j=\delta^k_j\delta^i_\ell$. One has
the
relations \begin{equation} E^k_\ell E^r_s=\delta^k_s E^r_\ell\ \ \mbox{and}\ \
\sum^N_{n=1} E^n_n=\bbbone \label{eq2} \end{equation} It follows from
(\ref{eq2})
that $M_N(\mathbb C)$ is a $\mathbb Z_N$-graded algebra if one equips it with
the
$\mathbb Z_N$-graduation defined by giving the degree $k-\ell$
$(\mbox{mod}(N))$
to $E^k_\ell$; $M_N(\mathbb C)=\displaystyle{\oplusinf_{p\in {\mathbb
Z}_N}}(M_N({\mathbb C}))^p$. Let $e=\lambda_1E^2_1+\dots
+\lambda_{N-1}E^N_{N-1}+\lambda_NE^1_N$, $(\lambda_i\in \mathbb C)$, be an
arbitrary element of degree 1, ($e\in (M_N({\mathbb C}))^1)$, i.e. \[ e= \left(
\begin{array}{c c c c c c} 0 &\lambda_1 & 0 & &\dots & 0\\ 0 & 0 & & &  & 0\\
\vdots &  & & & & \vdots\\ 0 &  && & & 0\\ 0 & & &&  & \lambda_{N-1}\\
\lambda_N
& 0 & 0 && \dots & 0 \end{array} \right) \] One has in view of (\ref{eq2})
\begin{equation} e^N=\lambda_1\lambda_2\dots \lambda_N\bbbone \label{eq3}
\end{equation} One defines a linear mapping of degree 1 of $M_N({\mathbb C})$
into itself by setting $d(A)=eA-q^aAe$ for $A\in (M_N({\mathbb C}))^a$. The
linear mapping $d$ satisfies the $q$-Leibniz rule \[ d(AB)=d(A) B+q^a Ad(B),\ \
\forall A\in (M_N({\mathbb C}))^a, \forall B\in M_N({\mathbb C}) \] Moreover
(\ref{eq3}) implies that $d^N=0$, (since $d^N=ad(e^N)$ as easily verified).
Thus
$M_N({\mathbb C})=\displaystyle{\oplusinf_p}(M_N({\mathbb C}))^p$ equipped with
$d$ satisfies the axioms of graded $q$-differential algebra except that it is
${\mathbb Z}_N$-graded instead of being ${\mathbb N}$-graded. However the
${\mathbb N}$-graded covering $p^\ast M_N({\mathbb C})$ equipped with
$p^\ast(d)$, (see in Section 1), is a graded $q$-differential algebra. The
algebra $(p^\ast M_N({\mathbb C}))^0=(M_N({\mathbb C}))^0$ of diagonal matrices
identifies with the algebra ${\mathbb C}^N$ of complex functions on a set with
$N$ elements and therefore the above graded $q$-differential algebra is a
$q$-differential calculus over the commutative algebra ${\mathbb C}^N$. Notice
that for $N=2$, $p^\ast M_2(\mathbb C)$ is an ordinary graded differential
algebra which is isomorphic to the universal differential envelope
$\Omega(\mathbb C^2)$ of the commutative algebra ${\mathbb C}^2$.

\subsection*{Example 2: Hochschild cochains}

Let $\cala$ be a unital associative ${\mathbb C}$-algebra and let ${\calm}$ be
a
bimodule over $\cala$. Recall that a $\calm$-valued Hochschild cochain of
degree
$n\in {\mathbb N}$ is a $n$-linear mapping of
$\underbrace{\cala\times\dots\times\cala}_n$ into $\calm$, (i.e. a linear
mapping
of $\otimes^n\cala$ into $\calm$). The vector space of $\calm$-valued
Hochschild cochains of degree $n$ is denoted by $C^n(\cala,\calm)$. The vector
space $C(\cala,\calm)=\displaystyle{\oplusinf_{n\in {\mathbb
N}}}C^n(\cala,\calm)$ of all $\calm$-valued Hochshild cochains is a ${\mathbb
N}$-graded vector space. If $\calm'$ is another bimodule over $\cala$, one
defines a graded bilinear mapping of $C(\cala,\calm)\times C(\cala,\calm')$
into
$C(\cala,\calm\displaystyle{\otimesinf_\cala}\calm'),\ (\alpha, \alpha')\mapsto
\alpha\cup\alpha'$, {\sl the cup product}, by setting for $\alpha\in
C^a(\cala,\calm)$ and $\alpha'\in C^{a'}(\cala,\calm')$ \[
(\alpha\cup\alpha')(x_1,\dots,x_{a+a'})=\alpha(x_1,\dots,x_a)
\displaystyle{\otimesinf_\cala}\alpha'(x_{a+1},\dots,x_{a+a'}),
\ \ \forall x_i\in \cala. \] The cup product if associative in the sense that
if
$\calm''$ is a third bi-module over $\cala$, one has for $\alpha\in
C(\cala,\calm),\alpha'\in C(\cala,\calm')$ and $\alpha''\in C(\cala,\calm'')$:
$(\alpha\cup\alpha')\cup \alpha''=\alpha\cup(\alpha'\cup\alpha'')$. By taking
$\calm=\calm'=\cala(=\calm'')$ and by making the identification
$\cala\displaystyle{\otimesinf_\cala}\cala=\cala$, one sees that, equipped with
the cup product, $C(\cala,\cala)$ is a unital $\mathbb N$-graded algebra with
$C^0(\cala,\cala)=\cala$. Let $q$ be a complex number with $q\not= 0$. One
defines a linear endormorphism $\delta_q$ of degree one of $C(\cala,\calm)$ by
setting for $\omega\in C^n(\cala,\calm)$, $\delta_1\omega=0$ and, for $q\not=
1$:
\begin{equation} \begin{array}{ll}
\delta_q(\omega)(x_0,\dots,x_n)=x_0\omega(x_1,\dots,x_n)& +\
\sum^n_{k=1}q^k\omega(x_0,\dots,x_{k-1}x_k,\dots,x_n)\\ \\ & -\
q^n\omega(x_0,\dots,x_{n-1})x_n \end{array} \label{eq4} \end{equation} $\forall
x_i\in \cala$. One verifies that $\delta^N_q=0$ whenever $q^N=1$ $(N\not=0)$
and
that, if $\beta\in C(\cala,\calm')$, one has:
$\delta_q(\omega\cup\beta)=\delta_q(\omega)\cup \beta + q^n\omega\cup
\delta_q(\beta)$. This implies in particular that $C(\cala,\cala)$ equipped
with
$\delta_q$ is a graded $q$-differential algebra and that it is a
$q$-differential
calculus over $\cala$. Notice that $\delta_{(-1)}$ is the usual Hochschild
coboundary $\delta$ so, when $q^N=1$ ($N\geq 2$), the
$H^{(k)}(C(\cala,\calm),\delta_q)$ defined as in Section 1 are $q$-analog of
Hochschild cohomology.

\subsection*{Example 3: $q$-differential dual of a product}

Let $\cala$ be an associative $\mathbb C$-algebra and let
$C(\cala)=\displaystyle{\oplusinf_{n\in {\mathbb N}}}C^n(\cala)$ be the graded
vector space of multilinear forms on $\cala$; i.e.
$C^n(\cala)=(\otimes^n\cala)^\ast$ is the ($\mathbb C$-) dual of
$\otimes^n\cala$
and $C^0(\cala)=\mathbb C$. By making the natural identifications
$C^n(\cala)\otimes C^m(\cala)\subset C^{n+m}(\cala)$ one sees that $C(\cala)$
is
canonically a $\mathbb N$-graded unital $\mathbb C$-algebra, (the product being
the tensor product over $\mathbb C$). By duality, the product ${\mathfrak
m}:\cala\otimes\cala\rightarrow \cala$ of $\cala$ gives a linear mapping
$\fm^\ast$ of $\cala^\ast$ into $(\cala\otimes \cala)^\ast$ i.e.
$\fm^\ast:C^1(\cala)\rightarrow C^2(\cala)$. For $q\in {\mathbb C}\backslash
\{0,1\}$, $\fm^\ast$ extends into a linear mapping
$\fm^\ast_q:C(\cala)\rightarrow C(\cala)$ which satisfies the graded
$q$-Leibniz
rule with \begin{equation} \fm^\ast_q(\omega)(x_0,\dots,x_n)=\sum^n_{k=1}
q^{k-1}\omega (x_0,\dots, x_{k-1} x_k,\dots,x_n) \label{eq5} \end{equation} for
$\omega\in C^n(\cala)$ and $x_i\in \cala$. It follows then from the
associativity
of the product of $\cala$ that one has $(\fm^\ast_q)^N=0$ whenever $q^N=1$,
$N\in
{\mathbb N}\backslash \{0\}$. Thus $C(\cala)$ equipped with $\fm^\ast_q$ is a
graded $q$-differential algebra. It is worth noticing here that the $\delta_q$
defined by (\ref{eq4}) on $C(\cala,\calm)$ in  Example 2 is up to a factor $q$
the $\fm^\ast_q$ defined by (\ref{eq5}), (i.e. ``the dual" of the product of
$\cala$) combined with a ``$q$-twisted bimodule action" on $\calm$. It should
also
be stressed that the results in Example 2 are true if $\cala$ is not unital
except that then $C(\cala,\cala)$ is also not unital.\\

Let us now drop the assumption that $\cala$ is associative, i.e. let $\cala$ be
a
complex vector space equipped with a bilinear product $(x,y)\mapsto xy$. Then
the
formula (\ref{eq5}) still defines a homogeneous linear mapping $\fm^\ast_q$ of
degree 1 of $C(\cala)$ into itself satisfying the $q$-Leibniz rule which
extends
the dual of the product, but now $q^N=1$ does not imply $(\fm^\ast_q)^N=0$. Let
$q$ be a $N$-th primitive root of the unity with $N\geq 2$. For $N=2$,
$(\fm^{\ast}_{(-1)})^2=0$ is equivalent to the associativity of the product of
$\cala$. For $N\geq 3$, $(\fm^\ast_q)^N=0$ is equivalent to a generalization of
degree $N+1$ of the associativity of the product of $\cala$ which is of the
form
$R_q(x_0\otimes x_1\otimes \dots \otimes x_N)=0$, $\forall x_i\in \cala$, where
$R_q$ is a linear mapping of $\otimes^{N+1}\cala$ into $\cala$. However, it
was remarked by Peter~W.~Michor \cite{PWM} that, if $\cala$ has a unit, then
the
relation $R_q=0$ implies the associativity of the product of $\cala$, i.e.
$R_{(-1)}=0$. Let us prove this fact. So let us assume that there is a
$\bbbone\in \cala$ such that $\bbbone x=x\bbbone=x,\ \forall x\in \cala$, and
let
$q$ be a $N$-th primitive root of the unity with $N\geq 3$. Then one has for
$x,y,z\in \cala$ and $\omega\in C^1(\cala)(=\cala^\ast)$,
$(\fm^\ast_q)^N\omega(x,y,\underbrace{\bbbone,\dots,\bbbone}_{N-2},z)=[N-2]_q
q^{N-2}\omega((x,y)z-x(yz))$. Since $\omega$ is arbitrary this shows that
$(\fm^\ast_q)^N=0$ implies the associativity of the product of $\cala$, (we
already know that the associativity of the product of $\cala$ implies
$(\fm^\ast_q)^N=0$ whenever $q^N=1$, $q\not= 1$ and $N\in {\mathbb N}\backslash
\{ 0\}$).

\section{The tensor algebra over $\cala$ of $\cala\otimes \cala$}

In this section $\cala$ is a unital associative $\mathbb C$-algebra. The tensor
product (over $\mathbb C$) $\cala\otimes\cala$ is in a natural way a bimodule
over $\cala$. The tensor algebra over $\cala$ of the bimodule
$\cala\otimes\cala$
will be denoted by $\mathfrak{T}(\cala)=\displaystyle{\oplusinf_{n\in{\mathbb
N}}}\mathfrak{T}^n(\cala)$. This is a unital graded algebra with
$\mathfrak{T}^n(\cala)=\otimes^{n+1}\cala$ and product defined by \[
(x_1\otimes\dots\otimes x_m)(y_1\otimes \dots\otimes
y_n)=x_1\otimes\dots\otimes
x_{m-1}\otimes x_m y_1\otimes y_2\otimes \dots \otimes y_n \ \mbox{for}\
x_i,y_j\in \cala. \]

In particular $\cala$ coincides with the subalgebra $\mathfrak{T}^0(\cala)$. As
a
tensor algebra over $\cala$, $\mathfrak{T}(\cala)$ satisfies  a universal
property. Here, since $\cala$ is unital, $\cala\otimes\cala$ is the free
bimodule
generated by $\tau=\bbbone \otimes\bbbone$. Hence $\mathfrak{T}(\cala)$ is the
$\mathbb N$-graded algebra generated by $\cala$ in degree 0 and by a free
generator $\tau$ of degree 1. In fact one has $x_0\otimes\dots\otimes
x_n=x_0\tau
x_1\dots \tau x_n$, $\forall x_i\in \cala$. Thus the graded algebra
$\mathfrak{T}(\cala)$ is also characterized by the following property.

\begin{lemma} Let $\mathfrak{A}=\oplus \mathfrak{A}^n$ be a unital $\mathbb
N$-graded $\mathbb C$-algebra, then for any homomorphism
$\varphi:\cala\rightarrow \mathfrak{A}^0$ of unital algebras and for any
$\alpha\in \mathfrak{A}^1$, there is a unique homorphism
$\mathfrak{T}_{\varphi,\alpha}:\mathfrak{T}(\cala)\rightarrow \mathfrak{A}$ of
graded algebras which extends $\varphi$, (i.e.
$\mathfrak{T}_{\varphi,\alpha}\restriction \cala=\varphi$), and is such that
$\mathfrak{T}_{\varphi,\alpha}(\tau)=\alpha$. \end{lemma} As an example of
application of this lemma, let us take $\mathfrak{A}=C(\cala,\cala)$, i.e. the
algebra of $\cala$-valued cochains of $\cala$ (see Example 2 of Section 2),
take
for $\varphi$ the identity mapping of $\cala$ onto itself considered as a
homomorphism of $\cala$ into $C^0(\cala,\cala)$ and take (again) for $\alpha$
the
identity mapping of $\cala$ onto itself considered as an element of
$C^1(\cala,\cala)$. Let $\Psi=\mathfrak{T}_{\varphi,
\alpha}:\mathfrak{T}(\cala)\rightarrow C(\cala,\cala)$ be the corresponding
graded-algebra homomorphism. This homomorphism which was considered in
\cite{Th.M} is given by \begin{equation} \Psi(x_0\otimes \dots \otimes
x_n)(y_1,\dots,y_n)=x_0y_1x_1\dots y_nx_n \label{eq6} \end{equation}

We now equip $\mathfrak{T}(\cala)$ with a structure of graded $q$-differential
algebra. Let $q$ be a complex number different from 0 and 1. One has the
following lemma \begin{lemma} There is a unique linear mapping
$d_q:\mathfrak{T}(\cala)\rightarrow \mathfrak{T}(\cala)$ homogeneous of degree
1
satisfying the $q$-Leibniz rule such that \[ d_q(x)=\bbbone\otimes x-x\otimes
\bbbone=\tau x - x\tau,\ \forall x\in \cala, \] and \[ d_q(\tau)=\tau^2,\
(\mbox{i.e.}\ d_q(\bbbone\otimes \bbbone)=\bbbone \otimes \bbbone \otimes
\bbbone). \] Moreover $d_q$ satisfies $d_q^N=0$ whenever $q^N=1$ for $N\geq 2$,
$N\in \mathbb N$. \end{lemma}

\noindent \underbar{Proof} It follows from the very structure of
$\mathfrak{T}(\cala)$ that for any derivation $D$ of $\cala$ into the bimodule
$\mathfrak{T}^1(\cala)=\cala\otimes\cala$ and for any $\mu\in
\mathfrak{T}^2(\cala)=\cala\otimes\cala\otimes\cala$, there is a unique
$D_q:\mathfrak{T}(\cala)\rightarrow \mathfrak{T}(\cala)$ satisfying the
$q$-Leibniz rule and such that $D_q(x)=D(x)$ for $x\in \cala$ and
$D_q(\tau)=\mu$. The first part of the lemma follows since $\mbox{ad}\ (\tau)$
is
a derivation of $\cala$ into $\mathfrak{T}^1(\cala)$. On the other hand, by
induction on $N\in \mathbb N$ with $N\geq 1$, one has for $x\in \cala$
\begin{equation} d^N_q(x)=[N!]_q\tau^{N-1}d_q(x)\ \ \mbox{and}\ \
d^N_q(\tau)=[N!]_q\tau^{N+1} \label{eq7} \end{equation} so the remaining part
of
the lemma follows from $[N]_q=0$ whenever $q^N=1$ for $N\in \mathbb N$ with
$N\geq 2$.$\square$\\

Thus $\mathfrak{T}(\cala)$ equipped with $d_q$ is a graded $q$-differential
algebra and, in fact a $q$-differential calculus over $\cala$. One verifies
that,
if $C(\cala,\cala)$ is equipped with the $\delta_q$ given by (\ref{eq4}), then
the above homomorphism $\Psi:\mathfrak{T}(\cala)\rightarrow C(\cala,\cala)$
given
by (\ref{eq6}) is an homomorphism of graded $q$-differential algebra, i.e. one
has $\Psi\circ d_q=\delta_q\circ\Psi$. This generalizes the result of
\cite{Th.M}
which is the case $q=-1$.

\subsection*{Remark 1}

There is another natural $q$-differential $d'_q$ on $\mathfrak{T}(\cala)$ which
is defined to be the unique linear mapping of $\mathfrak{T}(\cala)$ into itself
satisfying the $q$-Leibniz rule such that $d'_q(x)=d_q(x)=\tau x-x\tau$ for
$x\in\cala$ and $d'_q(\tau)=-q\tau^2$. One verifies that $d^{\prime N}_q=0$
whenever $q^N=1$ for $N\geq 2, \  N\in \mathbb{N}$. Correspondingly, there is
another $q$-differential $\delta'_q$ on $C(\cala,\cala)$ which, instead of
formula (\ref{eq4}), is given by \[ \begin{array}{ll}
\delta'_q(\omega)(x_0,\dots,x_n)=x_0\omega(x_1,\dots,x_n)& -\
\sum^n_{k=1}q^{k-1}\omega(x_0,\dots,x_{k-1}x_k,\dots,x_n)\\ \\ & -\
q^n\omega(x_0,\dots,x_{n-1})x_n \end{array} \] and is such that $\Psi\circ
d'_q=\delta'_q\circ\Psi$. The same formula gives, more generally, an
endomorphism
of $C(\cala,\calm)$ for any bimodule $\calm$ which has the same properties as
$\delta_q$. Notice that all these definitions coincide when $q=-1$.

\subsection*{Remark 2}

It is worth noticing here that both $q$-differentials $d_q$ and $d'_q$ on
$\mathfrak{T}(\cala)$ coincide on $\cala=\mathfrak{T}^0(\cala)$ with the
universal derivation $d:\cala\rightarrow \Omega^1(\cala)\subset
\mathfrak{T}^1(\cala)$, (see in Section 1). this is natural in view of the fact
that we shall represent the universal $q$-differential envelope of $\cala$ as a
$q$-differential subalgebra of $\mathfrak{T}(\cala)$. Now given an arbitrary
$\mu\in \mathfrak{T}^2(\cala)$ there is a unique $\tilde d$ on
$\mathfrak{T}(\cala)$ satisfying the $q$-Leibniz rule which extends the
universal
derivation $d$ and is such that $\tilde d(\tau)=\mu$. However in general one
does
not have $\tilde d^N=0$ when $q$ is a $N$-th primitive root of the unity. The
choices $\mu=\bbbone\otimes \bbbone \otimes \bbbone=\tau^2$ and $\mu=-q\bbbone
\otimes \bbbone \otimes \bbbone=-q\tau^2$ are the two choices for which this
generically holds.

\subsection*{Remark 3}

As a graded algebra, the universal differential envelope $\Omega(\cala)$ of
$\cala$ is a graded subalgebra of $\mathfrak{T}(\cala)$. The space
$\Omega^n(\cala)$ is the subspace of $\mathfrak{T}^n(\cala)$ which is
annihilated
by applying the multiplication $\fm$ of $\cala$ to two consecutive arguments.
On
the other hand, for $q=-1$, $d_{(-1)}$ is an ordinary differential on
$\mathfrak{T}(\cala)$ for which $\Omega(\cala)$ is stable. In fact,
$\Omega(\cala)$ {\sl is the smallest differential subalgebra of
$\mathfrak{T}(\cala)$ equipped with $d_{(-1)}$ which contains $\cala$.} We
shall
generalize this result by showing that the universal $q$-differential envelope
of
$\cala$ can be identified with the smallest $q$-differential subalgebra of
$\mathfrak{T}(\cala)$ equipped with $d_q$ which contains $\cala$ (i.e. the
$q$-differential subalgebra of $\mathfrak{T}(\cala)$ generated by $\cala$).

\section{Universal $q$-differential envelope}

In this section $\cala$ is a unital associative $\mathbb C$-algebra with unit
denoted by $\bbbone$ and $q\in \mathbb C\backslash \{0\}$ as before. If $q$ is
a
root of the unity, we define $N$ to be the smallest strictly positive integer
such that $q^N=1$, otherwise we set $N=\infty$. Let $d^k(\cala)$ for
$k\in\{1,2,\dots,N-1\}$ be $N-1$ copies of the vector space
$\cala/\mathbb{C}\bbbone$, $d^k:\cala\rightarrow d^k(\cala)$ being the
corresponding canonical projections. We extend $d:\cala\rightarrow d(\cala)$ as
a
linear mapping, again denoted by $d$, of $\cala\oplus d(\cala)\oplus
d^2(\cala)\oplus\dots\oplus d^{N-1}(\cala)$ into itself by defining
$d:d^k(\cala)\rightarrow d^{k+1}(\cala)$ to be the canonical isomorphism for
$k=1,2,\dots,N-2$ and by $d(d^{N-1}(\cala))=0$. The space $\cala\oplus
{\displaystyle{\oplus^{N-1}_{k=1}}}d^k(\cala)$ is equipped with a structure of
graded vector space by giving the degree 0 to the elements of $\cala$ and the
degree $k$ to the elements of $d^k(\cala)$ for $k=1,2,\dots,N-1$. The
endomorphism
$d$ is homogeneous of degree 1 and the graded subspace
$\cale={\displaystyle{\oplus^{N-1}_{k=1}}}d^k (\cala)$ is preserved by $d$.
Notice that the canonical projection $d^k:{\cal A}\rightarrow d^k({\cal A})$
coincides then with $\underbrace{d\circ\dots\circ d}_k:{\cal A}\rightarrow
d^k(\cal A)$, etc. so the notations are coherent. Let $T(\cale)$ be the tensor
algebra over the graded vector space
$\cale={\displaystyle{\oplus^{N-1}_{k=1}}}d^k({\cala})$. On $T({\cale})$ there
is a unique graduation compatible with the graduation of ${\cale}$ such that it
is a graded algebra and on this graded algebra there is a unique extension,
again
denoted by $d$, of the endormorphism $d$ of ${\cale}$ which satisfies the
$q$-Leibniz rule. Namely one has for $x_i\in {\cala}$ and
$k_i\in\{1,\dots,N-1\}$
\[ \begin{array}{l} \partial(d^{k_1}(x_1)\dots d^{k_n}(x_n)) = k_1+\dots
+k_n,\\
\\ d(d^{k_1}(x_1)\dots d^{k_n}(x_n)) =\\ \\
\displaystyle{\sum^n_{i=1}}q^{k_1+\dots +k_{i-1}}d^{k_1}(x_1)\dots
d^{k_{i-1}}(x_{i+1})d^{k_i+1}(x_i)d^{k_{i+1}}(x_{i+1})\dots d^{k_n}(x_n)
\end{array} \] where $\partial$ denotes the degree and the product is the
tensor
product. Formula (1) is satisfied therefore, for $N<\infty$,
$d^N=\underbrace{d\circ\dots\circ d}_{N}$ is a derivation which vanishes on
$T(\cale)$ since it vanishes on $\cale$. Thus $T(\cale)$ is a graded
$q$-differential algebra.\\

Let $\Omega_q(\cala)$ be defined by $\Omega_q(\cala)=\cala\otimes T(\cale)$.
The
space $\Omega_q(\cala)$ is a graded vector space with graduation given by
$\partial(x\otimes t)=\partial(t)$ for $x\in \cala$ and $t\in T(\cale)$. It is
also canonically a left $\cala$-module and a graded right $T(\cale)$-module for
the above graduation. One extends all the previous definitions of $d$ to
$\Omega_q(\cala)$ by setting $d(x\otimes t)=\bbbone \otimes d(x) t+x\otimes
d(t)=d(x)t+xd(t)$ for $x\in \cala,\ t\in T(\cale)$ and where in the last
equality
$\bbbone \otimes T(\cale)$ and $T(\cale)$ are identified. The endomorphism $d$
satisfies in $\Omega_q(\cala)$ $d((x\otimes t)t')=d(x\otimes
t)t'+q^{\partial(x\otimes t)}(x\otimes t)d(t')$ for $x\otimes t\in
\Omega_q(\cala)$ and $t'\in T(\cale)\ (\subset\Omega_q(\cala))$. Identifying
$\cala$ with $\cala\otimes\bbbone\subset\Omega_q(\cala)$, one has the
following.
\begin{lemma} There is a unique associative product on $\Omega_q(\cala)$ which
extends its structure of $(\cala,T(\cale))$-bimodule, for which
$\Omega_q(\cala)$
is a graded algebra and for which $d$ satisfies the $q$-Leibniz rule. Then
$\Omega_q(\cala)$ equipped with $d$ is a graded $q$-differential algebra, (i.e.
$d^N=0$ for $N<\infty$). \end{lemma} \noindent\underbar{Proof.} What is needed
is
a product on the right by elements of $\cala$. If the $q$-Leibniz is satisfied
for $d$, one must have (formula (\ref{eq1})) \[d^n(x)b=d^n(xb)-\sum^n_{p=1}
\left[\begin{array}{c}n\\p\end{array}\right]_q d^{n-p}(x)d^p(b)\in
\Omega_q(\cala),\ \ \forall x, b\in \cala\] and $\forall n\geq 1$, from which
the
uniqueness of the product follows if it exists. Define $(yd^n(x))b$ for
$y,x,b\in
\cala$ and $n\geq 1$ by the above formula i.e.
\[(yd^n(x))b=yd^n(xb)-\sum^n_{p=1}
\left[\begin{array}{c}n\\ {p}\end{array}\right]_q yd^{n-p}(x)d^p(b).\] By
definition, one has $(yd^n(x))b=y(d^n(x)b)$. On the other hand, it follows from
the properties of the $q$-binomial coefficients that one has for $y,x,b,c\in
\cala$ and $n\geq 1$ $(yd^n(x))(bc)=((yd^n(x))b)c$ and that therefore the
product
extends uniquely into an associative one by setting $(yt)(y't')=((yt)y')t'$ for
$y,y'\in \cala$ and $t,t'\in T(\cale)$. The fact that $d$ satisfies the
$q$-Leibniz rule follows from
$d^n(x)b=d^k(d^{n-k}(x)b)-\displaystyle{\sum^k_{p=1}}q^{(n-k)p}
\left[\begin{array}{c}
k\\p\end{array}\right]_q d^{n-p}(x)d^p(b)$, for $x,b\in \cala$ and $n\geq k\geq
1$.  Furthermore, for $N<\infty$, one has $d^N=0$ since $d^N$ vanishes on the
generators. Thus $\Omega_q(\cala)$ is a graded $q$-differential
algebra.$\square$

\begin{theorem} Let $\fraca=\displaystyle{\oplusinf_{n\in \mathbb N}}\fraca^n$
be
a graded $q$-differential algebra and let $\varphi:\cala\rightarrow \fraca^0$
be
a homomorphism of unital algebras. Then there is a unique homomorphism
$\bar\varphi:\Omega_q(\cala)\rightarrow \fraca$ of graded $q$-differential
algebras which induces $\varphi$. \end{theorem}

\noindent\underbar{Proof.} In any graded $q$-differential algebra one has
$d(\bbbone)=0$, therefore one defines a linear mapping
$\varphi_0:\cale\rightarrow \fraca$ by setting $\varphi_0(d^kx)=d^k\varphi(x)$
for $x\in \cala$ and $k\geq 1$. By the universal property of the tensor algebra
$\varphi_0$ extends uniquely into an algebra homomorphism
$\varphi_1:T(\cale)\rightarrow \fraca$. The homomorphism $\varphi_1$ is
obviously
a homomorphism of graded algebras satisfying $\varphi_1\circ d=d\circ
\varphi_1$
so it is a homomorphism of graded $q$-differential algebras. Define the linear
mapping $\bar\varphi:\Omega_q(\fraca)\rightarrow\fraca$ by
$\bar\varphi(xt)=\varphi(x)\varphi_1(t)$ for $x\in\cala$ and $t\in T(\cale)$.
One
has $\bar\varphi(x)=\varphi(x)$ for $x\in\cala$, $\bar\varphi((xt)t')=\bar
\varphi(xt)\bar \varphi(t')$ for $t,t'\in T(\cale)$, $\bar\varphi\circ d=d\circ
\bar\varphi$ and $\bar\varphi$ is unique under these conditions. It follows
that
$\bar\varphi$ is in fact a homomorphism of graded $q$-differential algebras
which
is unique under the condition that
$\bar\varphi\restriction\cala=\varphi$.$\square$\\

The graded $q$-differential algebra $\Omega_q(\cala)$ is characterized uniquely
up to an isomorphism by the universal property stated in Theorem 1, this is why
$\Omega_q(\cala)$ will be called {\sl the universal $q$-differential envelope
of
$\cala$} or {\sl the universal $q$-differential calculus over $\cala$}.
\begin{proposition} The canonical homomorphism
$\overline{\mbox{Id}}:\Omega_q(\cala)\rightarrow \mathfrak{T}(\cala)$ induced
by
the identity mapping of $\cala$ onto itself (as in Theorem 1) is injective
$(q\not=0$ and $q\not=1)$. \end{proposition}

\noindent\underbar{Proof} In $\mathfrak{T}(\cala)$ one has (\ref{eq7}) \[
d^k_q(x)=[k!]_q(\bbbone^{\otimes^k}\otimes x-\bbbone^{\otimes^{k-1}}\otimes
x\otimes \bbbone)\ \ \mbox{for}\  \ k\in\{1,2,\dots,N-1\},\ x\in\cala. \] This
implies that $\overline{\mbox{Id}}$ induces an isomorphism of $T(\cale)$ onto
the
subalgebra of $\mathfrak{T}(\cala)$ generated by the $d^k_q(x)$ for
$k\in\{1,\dots,N-1\}$ and $x\in\cala$. The remaining follows from the fact that
the left $\cala$-submodule of $\mathfrak{T}(\cala)$ generated by the image of
$T(\cale)$, is freely generated, i.e. is isomorphic to $\cala\otimes
T(\cale)$.$\square$\\

Thus, one can identify $\Omega_q(\cala)$ with the $q$-differential subalgebra
of
$\mathfrak{T}(\cala)$ generated by $\cala$. This generalizes the standard
representation of the usual universal differential envelope of $\cala$,
\cite{Kar}, which is the case $q=-1$.\\

There is another approach of the construction of $\Omega_q(\cala)$ as
$q$-differential subalgebra of $\ft(\cala)$ which we now sketch. This approach
is
based on the universal Hochschild cocycles \cite{CKMV}, \cite{CuQ}. Recall that
a
derivation $X$ of $\cala$ into a bimodule $\calm$ is a $\calm$-valued
Hochschild
1-cocycle. If $X_i:\cala\rightarrow \calm_i$ are $n$ derivations of $\cala$
into
bimodules $\calm_i$, then their cup product $X_1\cup \cdots \cup
X_n:\otimes^n\cala\rightarrow \calm_1\otimes_{\cala}\cdots
\otimes_{\cala}\calm_n$ is a $\calm_1\otimes_\cala
\cdots\otimes_\cala\calm_n$-valued Hochschild $n$-cocycle. This cocycle is
normalized in the sense that it vanishes whenever one of its arguments is the
unit $\bbbone$ of $\cala$. Consider in particular the universal derivation
$d:\cala\rightarrow \Omega^1(\cala)$, (see in Section 1). By taking the cup
product $n$ times with itself of $d$, one obtains a normalized $n$-cocycle
$d^{\cup^n}:\otimes^n\cala\rightarrow \Omega^n(\cala)$ which is defined by
$d^{\cup^n}(x_1,\cdots ,x_n)=d(x_1)\cdots d(x_n)$. It turns out that this
normalized $n$-cocycle is universal, \cite{CKMV}, \cite{CuQ}, in the sense that
{\sl for any normalized Hochschild $n$-cocycle $c:\otimes^n\cala\rightarrow
\calm$ of $\cala$ into a bimodule $\calm$, there is a unique bimodule
homomorphism $i_c$ of $\Omega^n(\cala)$ into $\calm$ such that} $c=i_c\circ
d^{\cup^n}$. Furthermore one can characterize the triviality of $c$ in terms of
the homomorphism $i_c$, \cite{CKMV}. For that one notices that one has the
inclusions of bimodules $\Omega^n(\cala)\subset \cala\otimes
\Omega^{n-1}(\cala)\subset \ft^n(\cala)$ for $n\geq 1$. More precisely one has
an
exact sequence $(n\geq 1)$ \[ 0\rightarrow \Omega^n(\cala)\fc
\cala\otimes\Omega^{n-1}(\cala) \displaystyle{\buildrel \fm\over
\longrightarrow}\Omega^{n-1}(\cala)\rightarrow 0. \] where $\fm$ is the
multiplication of $\Omega(\cala)$. The cocycle $d^{\cup^n}$ is non-trivial,
however it is trivial if it is considered as a $\cala\otimes
\Omega^{n-1}(\cala)$-valued cocycle because one has there~: \[
\begin{array}{ll}
dx_1\cdots dx_n=-(x_1\otimes dx_2\cdots dx_n & +\sum^{n-1}_{k=1}(-1)^k\bbbone
\otimes dx_1\cdots d(x_kx_{k+1})\cdots dx_n\\ \\ &+(-1)^n\bbbone \otimes
(dx_1\cdots dx_{n+1})x_n) \end{array} \] i.e. $d^{\cup^n}=\delta(-\bbbone
\otimes
d^{\cup^{n-1}})$ in $\cala\otimes \Omega^{n-1}(\cala)$, where $\delta$ is the
Hochschild coboundary $(\delta=\delta_{(-1)})$. Therefore, if the
$\calm$-valued
normalized $n$-cocycle $c$ is such that $i_c$ is the restriction to
$\Omega^n(\cala)$ of a bimodule homomorphism $\varphi:\cala\otimes
\Omega^{n-1}(\cala)\rightarrow \calm$, then it is trivial because one has
$c=\delta(\varphi(-\bbbone\otimes d^{\cup^{n-1}})$). Conversely, if
$c=\delta(c')$ then by setting $c'=\varphi(- \bbbone \otimes d^{\cup^{n-1}})$
one
defines an extension $\varphi$ of $i_c$ to $\cala\otimes \Omega^{n-1}(\cala)$.
Let us apply this to the construction of $\Omega_q(\cala)$. So let $q$ be a
complex number different from 0 and 1 and let $\fraca=\oplus \fraca^n$ be an
arbitrary graded $q$-differential algebra with $\fraca^0=\cala$. As already
stressed, the $q$-differential $d_{\fraca}$ of $\fraca$ induces a derivation of
$\cala$ into $\fraca^1$ so one must take $\Omega^1_q(\cala)=\Omega^1(\cala)$
and
the $q$-differential of $\Omega_q(\cala)$ must induce the universal derivation
$d:\cala\rightarrow \Omega^1(\cala)$. Then the normalized 2-cocycle
$d_\fraca\cup
d_\fraca$ induces a unique bimodule homomorphism $i_2=i_{d_\fraca\cup
d_\fraca}$
of $\Omega^2(\cala)$ into $\fraca^2$ so $\Omega^2(\cala)\subset
\Omega^2_q(\cala)$. However the $q$-Leibniz rules implies
$d^2_\fraca(xy)=xd^2_\fraca(y)+d^2_\fraca(x)y+(1+q)d_\fraca(x)d_\fraca(y)$. So
if
$q\not= -1$ then\linebreak[4] $d_\fraca\cup
d_\fraca=\delta\left(-\frac{1}{1+q}d^2_\fraca\right)=\delta\left(-\frac{1}{[2]_q}d^2_\fraca\right)$.
Therefore (if $q\not=-1$), by the above discussion, $i_2$ has a unique
extension
as a bimodule homomorphism\linebreak[4]
$\varphi:\cala\otimes\Omega^1(\cala)\rightarrow \fraca^2$ such that
$d^2_\fraca(x)=\varphi([2]_q\bbbone \otimes d(x))$, $\forall x\in \cala$. It
follows that, if $q\not=-1$, one must take $\Omega^2_q(\cala)=\cala\otimes
\Omega^1(\cala)\ (\subset \ft^2(\cala))$ and $d^2(x)=[2]_q\bbbone \otimes
d(x)=[2]_q\tau d(x)$ which is, in view of (7), the formula induced by the
$q$-differential $d_q$ of $\ft(\cala)$. Although it becomes a little
cumbersome,
one can continue the construction of $\Omega_q(\cala)$ as $q$-differential
subalgebra of $\ft(\cala)$ along this line, (by using the formula (1) and the
universal cocycles, etc).

\section{Conclusion}

In this paper we have generalized several constructions of ordinary
differential
algebra to $q$-differential algebra. When $q$ is a primitive $N$-th root of the
unity, (e.g. $q=\exp \left(\frac{2\pi i}{N}\right)$), with $N\geq 2$, it is
natural to ask what is the generalized cohomology $H^{(p),n}$ $(p=1,\dots,
N-1$,
$n\in \mathbb N)$ of the various graded $q$-differential algebras introduced
here. The computation of these generalized cohomologies will be described in a
separate paper \cite{M.D-V}, we just give here the results. For the graded
$q$-differential algebras $(C(\cala),\mathfrak m^\ast_q), (\mathfrak T(\cala),
d_q)$
and $\Omega_q(\cala)$ of Example 3, of Section 3 and of Section 4, these
generalized cohomologies are trivial as expected, i.e.  one has $H^{(p),n}=0$
for $n\geq 1$
and $H^{(p),0}=\mathbb C$, $p\in \{1,2,\dots, N-1\}$. For the case of the
generalized Hochschild cohomology i.e. of $(C(\cala,\calm),\ \delta_q)$ of
example 2 the result is the following: If $\cala$ {\sl is unital}, then one has
$H^{(p),Nk}=H^{2k}$ and $H^{(p),N(k+1)-p}=H^{2(k+1)-1}$ for
$p\in\{1,\dots,N-1\}$ and $k\in \mathbb N$, where $H^n$ denotes the usual
Hochschild cohomology, and
$H^{(p),r}=0$ otherwise i.e. if $r\not=0\ \mbox{mod}(N)$  and $r+p\not=0\
\mbox{mod}(N)$. Thus, for unital algebras, the information contained in the
generalized Hochschild cohomology is the same as the one of ordinary Hochschild
cohomology.

\section*{Acknowledgements}

One of the authors (M.D-V) is greatly indebted to Peter W. Michor for numerous
stimulating discussions and to Max Karoubi for his kind interest and
communication of reference \cite{Kapr}.

\end{document}